\documentclass{amsart}
\usepackage{amsfonts}

\newtheorem{theorem}{Theorem}
\theoremstyle{plain}

\newtheorem{corollary}{Corollary}

\newtheorem{proposition}{Proposition}
\newtheorem{remark}{Remark}

\numberwithin{equation}{section}
\input{tcilatex}

\begin{document}
\title[Schwarz's Inequality]{A Counterpart of Schwarz's Inequality in Inner
Product Spaces}
\author{Sever S. Dragomir}
\address{School of Computer Science and Mathematics\\
Victoria University of Technology\\
PO Box 14428, MCMC 8001\\
Victoria, Australia.}
\email{sever@matilda.vu.edu.au}
\date{May 21, 2003}
\urladdr{http://rgmia.vu.edu.au/SSDragomirWeb.html}

\begin{abstract}
A new counterpart of Schwarz's inequality in inner product spaces and
applications for isotonic functionals, integrals and sequences are provided.
\end{abstract}

\maketitle

\section{Introduction}

Let $\overline{a}=\left( a_{1},\dots ,a_{n}\right) $ and $\overline{b}%
=\left( b_{1},\dots ,b_{n}\right) $ be two positive $n-$tuples with%
\begin{equation}
0<m_{1}\leq a_{i}\leq M_{1}<\infty \text{ and }0<m_{2}\leq b_{i}\leq
M_{2}<\infty ;  \label{1.1a}
\end{equation}%
for each \ $i\in \left\{ 1,\dots ,n\right\} ,$ and some constants $%
m_{1},m_{2},M_{1},M_{2}.$

The following counterparts of the Cauchy-Bunyakowsy-Schwarz inequality are
valid:

\begin{enumerate}
\item P\'{o}lya-Szeg\"{o}'s inequality \cite{6b}%
\begin{equation}
1\leq \frac{\sum_{k=1}^{n}a_{k}^{2}\sum_{k=1}^{n}b_{k}^{2}}{\left(
\sum_{k=1}^{n}a_{k}b_{k}\right) ^{2}}\leq \frac{1}{4}\left( \sqrt{\frac{%
M_{1}M_{2}}{m_{1}m_{2}}}+\sqrt{\frac{m_{1}m_{2}}{M_{1}M_{2}}}\right) ^{2};
\label{1.2}
\end{equation}

\item Shisha-Mond's inequality \cite{7b}%
\begin{equation}
\frac{\sum_{k=1}^{n}a_{k}^{2}}{\sum_{k=1}^{n}a_{k}b_{k}}-\frac{%
\sum_{k=1}^{n}a_{k}b_{k}}{\sum_{k=1}^{n}b_{k}^{2}}\leq \left[ \left( \frac{%
M_{1}}{m_{2}}\right) ^{\frac{1}{2}}-\left( \frac{m_{1}}{M_{2}}\right) ^{%
\frac{1}{2}}\right] ^{2};  \label{1.3}
\end{equation}

\item Ozeki's inequality \cite{5b}%
\begin{equation}
\sum_{k=1}^{n}a_{k}^{2}\sum_{k=1}^{n}b_{k}^{2}-\left(
\sum_{k=1}^{n}a_{k}b_{k}\right) ^{2}\leq \frac{n^{2}}{4}\left(
M_{1}M_{2}-m_{1}m_{2}\right) ^{2};  \label{1.4}
\end{equation}

\item Diaz-Metcalf's inequality \cite{1b}%
\begin{equation}
\sum_{k=1}^{n}b_{k}^{2}+\frac{m_{2}M_{2}}{m_{1}M_{1}}\sum_{k=1}^{n}a_{k}^{2}%
\leq \left( \frac{M_{2}}{m_{1}}+\frac{m_{2}}{M_{1}}\right)
\sum_{k=1}^{n}a_{k}b_{k}.  \label{1.5}
\end{equation}
\end{enumerate}

If $\overline{\mathbf{w}}=\left( w_{1},\dots ,w_{n}\right) $ is a positive
sequence, then the following weighted inequalities also hold:

\begin{enumerate}
\item Cassel's inequality \cite{8b}. If the positive real sequences $%
\overline{a}=\left( a_{1},\dots ,a_{n}\right) $ and $\overline{b}=\left(
b_{1},\dots ,b_{n}\right) $ satisfy the condition 
\begin{equation}
0<m\leq \frac{a_{k}}{b_{k}}\leq M<\infty \text{ for each }k\in \left\{
1,...,n\right\}  \label{1.5a}
\end{equation}%
then%
\begin{equation*}
\frac{\left( \sum_{k=1}^{n}w_{k}a_{k}^{2}\right) \left(
\sum_{k=1}^{n}w_{k}b_{k}^{2}\right) }{\left(
\sum_{k=1}^{n}w_{k}a_{k}b_{k}\right) ^{2}}\leq \frac{\left( M+m\right) ^{2}}{%
4mM};
\end{equation*}

\item Greub-Reinboldt's inequality \cite{2b}%
\begin{equation}
\left( \sum_{k=1}^{n}w_{k}a_{k}^{2}\right) \left(
\sum_{k=1}^{n}w_{k}b_{k}^{2}\right) \leq \frac{\left(
M_{1}M_{2}+m_{1}m_{2}\right) ^{2}}{4m_{1}m_{2}M_{1}M_{2}}\left(
\sum_{k=1}^{n}w_{k}a_{k}b_{k}\right) ^{2}.  \label{1.6}
\end{equation}%
provided $\overline{a}=\left( a_{1},\dots ,a_{n}\right) $ and $\overline{b}%
=\left( b_{1},\dots ,b_{n}\right) $ satisfy the condition $\left( \ref{1.1a}%
\right) .$

\item Generalised Diaz-Metcalf inequality \cite{1b}, see also \cite[p. 123]%
{4b}. If $u,v\in \left[ 0,1\right] $ and $v\leq u,$ $u+v=1$ and $\left( \ref%
{1.5a}\right) $ holds, then one has the inequality%
\begin{equation}
u\sum_{k=1}^{n}w_{k}b_{k}^{2}+vMm\sum_{k=1}^{n}w_{k}a_{k}^{2}\leq \left(
vm+uM\right) \sum_{k=1}^{n}w_{k}a_{k}b_{k}.  \label{1.8}
\end{equation}

\item Klamkin-McLenaghan's inequality \cite{3b}. If $\overline{a},\overline{b%
}$ satisfy ($\ref{1.5a}$)$,$ then%
\begin{multline}
\left( \sum_{i=1}^{n}w_{i}a_{i}^{2}\right) \left(
\sum_{i=1}^{n}w_{i}b_{i}^{2}\right) -\left(
\sum_{i=1}^{n}w_{i}a_{i}b_{i}\right) ^{2}  \label{1.9} \\
\leq \left( M^{\frac{1}{2}}-m^{\frac{1}{2}}\right)
^{2}\sum_{i=1}^{n}w_{i}a_{i}b_{i}\sum_{i=1}^{n}w_{i}a_{i}^{2}.
\end{multline}
\end{enumerate}

For other results providing counterpart inequalities, see the recent
monograph on line \cite{1bab}.

In this paper we point out a new counterpart of Schwarz's inequality in real
or complex inner product spaces. Particular cases for isotonic linear
functionals, integrals and sequences are also provided.

\section{An Inequality in Inner Product Spaces}

The following reverse of Schwarz's inequality in inner product spaces holds.

\begin{theorem}
\label{t2.1}Let $A,a\in \mathbb{K}$ $\left( \mathbb{K}=\mathbb{C},\mathbb{R}%
\right) $ and $x,y\in H.$ If%
\begin{equation}
\func{Re}\left\langle Ay-x,x-ay\right\rangle \geq 0,  \label{2.1}
\end{equation}%
or, equivalently,%
\begin{equation}
\left\Vert x-\frac{a+A}{2}\cdot y\right\Vert \leq \frac{1}{2}\left\vert
A-a\right\vert \left\Vert y\right\Vert ,  \label{2.1.a}
\end{equation}%
holds, then one has the inequality%
\begin{equation}
0\leq \left\Vert x\right\Vert ^{2}\left\Vert y\right\Vert ^{2}-\left\vert
\left\langle x,y\right\rangle \right\vert ^{2}\leq \frac{1}{4}\left\vert
A-a\right\vert ^{2}\left\Vert y\right\Vert ^{4}.  \label{2.2}
\end{equation}%
The constant $\frac{1}{4}$ is sharp in (\ref{2.2}).
\end{theorem}

\begin{proof}
The equivalence between $\left( \ref{2.1}\right) $ and $\left( \ref{2.1.a}%
\right) $ can be easily proved, see for example \cite{1ba}.

Let us define%
\begin{equation*}
I_{1}:=\func{Re}\left[ \left( A\left\Vert y\right\Vert ^{2}-\left\vert
\left\langle x,y\right\rangle \right\vert \right) \left( \overline{%
\left\langle x,y\right\rangle }-\overline{a}\left\Vert y\right\Vert
^{2}\right) \right]
\end{equation*}%
and 
\begin{equation*}
I_{2}:=\left\Vert y\right\Vert ^{2}\func{Re}\left\langle
Ay-x,x-ay\right\rangle .
\end{equation*}%
Then%
\begin{equation*}
I_{1}=\left\Vert y\right\Vert ^{2}\func{Re}\left[ A\overline{\left\langle
x,y\right\rangle }+\overline{a}\left\langle x,y\right\rangle \right]
-\left\vert \left\langle x,y\right\rangle \right\vert ^{2}-\left\Vert
y\right\Vert ^{4}\func{Re}\left( A\overline{a}\right)
\end{equation*}%
and%
\begin{equation*}
I_{2}=\left\Vert y\right\Vert ^{2}\func{Re}\left[ A\overline{\left\langle
x,y\right\rangle }+\overline{a}\left\langle x,y\right\rangle \right]
-\left\Vert x\right\Vert ^{2}\left\Vert y\right\Vert ^{2}-\left\Vert
y\right\Vert ^{4}\func{Re}\left( A\overline{a}\right) ,
\end{equation*}%
giving%
\begin{equation}
I_{1}-I_{2}=\left\Vert x\right\Vert ^{2}\left\Vert y\right\Vert
^{2}-\left\vert \left\langle x,y\right\rangle \right\vert ^{2};  \label{2.3}
\end{equation}%
for any $x,y\in H$ and $a,A\in \mathbb{K}$, which is an interesting equality
in itself as well.

If (\ref{2.1}) holds, then $I_{2}\geq 0$ and thus%
\begin{equation}
\left\Vert x\right\Vert ^{2}\left\Vert y\right\Vert ^{2}-\left\vert
\left\langle x,y\right\rangle \right\vert ^{2}\leq \func{Re}\left[ \left(
A\left\Vert y\right\Vert ^{2}-\left\vert \left\langle x,y\right\rangle
\right\vert \right) \left( \overline{\left\langle x,y\right\rangle }-%
\overline{a}\left\Vert y\right\Vert ^{2}\right) \right] .  \label{2.4}
\end{equation}%
If we use the elementary inequality for $u,v\in \mathbb{K}$ $\left( \mathbb{K%
}=\mathbb{C},\mathbb{R}\right) $%
\begin{equation}
\func{Re}\left[ u\overline{v}\right] \leq \frac{1}{4}\left\vert
u+v\right\vert ^{2},  \label{2.5}
\end{equation}%
then we have for 
\begin{equation*}
u:=A\left\Vert y\right\Vert ^{2}-\left\langle x,y\right\rangle ,\ \
v:=\left\langle x,y\right\rangle -a\left\Vert y\right\Vert ^{2}
\end{equation*}%
that%
\begin{equation}
\func{Re}\left[ \left( A\left\Vert y\right\Vert ^{2}-\left\vert \left\langle
x,y\right\rangle \right\vert \right) \left( \overline{\left\langle
x,y\right\rangle }-\overline{a}\left\Vert y\right\Vert ^{2}\right) \right]
^{2}\leq \frac{1}{4}\left\vert A-a\right\vert ^{2}\left\Vert y\right\Vert
^{4}.  \label{2.6}
\end{equation}%
Making use of the inequalities $\left( \ref{2.4}\right) $ and $\left( \ref%
{2.6}\right) ,$ we deduce $\left( \ref{2.2}\right) $.

Now, assume that (\ref{2.2}) holds with a constant $C>0,$ i.e.,%
\begin{equation}
\left\Vert x\right\Vert ^{2}\left\Vert y\right\Vert ^{2}-\left\vert
\left\langle x,y\right\rangle \right\vert ^{2}\leq C\left\vert
A-a\right\vert ^{2}\left\Vert y\right\Vert ^{4},  \label{2.61}
\end{equation}%
where $x,y,a,A$ satisfy (\ref{2.1}).

Consider $y\in H,$ $\left\Vert y\right\Vert =1,$ $a\neq A$ and $m\in H,$ $%
\left\Vert m\right\Vert =1$ with $m\perp y.$ Define%
\begin{equation*}
x:=\frac{A+a}{2}y+\frac{A-a}{2}m.
\end{equation*}%
Then%
\begin{equation*}
\left\langle Ay-x,x-ay\right\rangle =\left\vert \frac{A-a}{2}\right\vert
^{2}\left\langle y-m,y+m\right\rangle =0
\end{equation*}%
and thus the condition (\ref{2.1}) is fulfilled. From (\ref{2.61}) we deduce%
\begin{equation}
\left\Vert \frac{A+a}{2}y+\frac{A-a}{2}m\right\Vert ^{2}-\left\vert
\left\langle \frac{A+a}{2}y+\frac{A-a}{2}m,y\right\rangle \right\vert
^{2}\leq C\left\vert A-a\right\vert ^{2},  \label{2.7}
\end{equation}%
and since%
\begin{equation*}
\left\Vert \frac{A+a}{2}y+\frac{A-a}{2}m\right\Vert ^{2}=\left\vert \frac{A+a%
}{2}\right\vert ^{2}+\left\vert \frac{A-a}{2}\right\vert ^{2}
\end{equation*}%
and%
\begin{equation*}
\left\vert \left\langle \frac{A+a}{2}y+\frac{A-a}{2}m,y\right\rangle
\right\vert ^{2}=\left\vert \frac{A+a}{2}\right\vert ^{2}
\end{equation*}%
then by (\ref{2.7}) we obtain%
\begin{equation*}
\frac{\left\vert A-a\right\vert ^{2}}{4}\leq C\left\vert A-a\right\vert ^{2},
\end{equation*}%
giving $C\geq \frac{1}{4},$ and the theorem is completely proved.
\end{proof}

\section{Applications for Isotonic Linear Functionals}

Let $F\left( T\right) $ be an algebra of real functions defined on $T$ and $%
L $ a subclass of $F\left( T\right) $ satisfying the conditions:

\begin{enumerate}
\item[(i)] $f,g\in L$ implies $f+g\in L;$

\item[(ii)] $f\in L,$ $\in \mathbb{R}$ implies $\alpha f\in L.$
\end{enumerate}

A functional $A$ defined on $L$ is an \textit{isotonic linear functional} on 
$L$ provided that

\begin{enumerate}
\item[(a)] $A\left( \alpha f+\beta g\right) =\alpha A\left( f\right) +\beta
A\left( g\right) $ for all $\alpha ,\beta \in \mathbb{R}$ and $f,g\in L;$

\item[(aa)] $f\geq g,$ that is, $f\left( t\right) \geq g\left( t\right) $
for all $t\in T,$ implies $A\left( f\right) \geq A\left( g\right) .$
\end{enumerate}

The functional $A$ is \textit{normalised} on $L,$ provided that $\mathbf{1}%
\in L,$ i.e., $\mathbf{1}\left( t\right) =1$ for all $t\in T,$ implies $%
A\left( \mathbf{1}\right) =1.$

Usual examples of isotonic linear functionals are integrals, sums, etc.

Now, suppose that $h\in F\left( T\right) ,$ $h\geq 0$ is given and satisfies
the properties that $fgh\in L,$ $fh\in L,$ $gh\in L$ for all $f,g\in L.$ For
a given isotonic linear functional $A:L\rightarrow \mathbb{R}$ with $A\left(
h\right) >0,$ define the mapping $\left( \cdot ,\cdot \right) _{A,h}:L\times
L\rightarrow \mathbb{R}$ by%
\begin{equation}
\left( f,g\right) _{A,h}:=\frac{A\left( fgh\right) }{A\left( h\right) }.
\label{3.2}
\end{equation}

This functional satisfies the following properties:

\begin{enumerate}
\item[(s)] $\left( f,f\right) _{A,h}\geq 0$ for all $f\in L;$

\item[(ss)] $\left( \alpha f+\beta g,k\right) _{A,h}=\alpha \left(
f,k\right) _{A,h}+\beta \left( g,k\right) _{A,h}$ for all $f,g,k\in L$ and $%
\alpha ,\beta \in \mathbb{R}$;

\item[(sss)] $\left( f,g\right) _{A,h}=\left( g,f\right) _{A,h}$ for all $%
f,g\in L.$
\end{enumerate}

The following reverse of Schwarz's inequality for positive linear
functionals holds.

\begin{proposition}
\label{p3.1}Let $f,g,h\in F\left( T\right) $ be such that $fgh\in L,$ $%
f^{2}h\in L,$ $g^{2}h\in L.$ If $m,M$ are real numbers such that%
\begin{equation}
mg\leq f\leq Mg\text{ on }F\left( T\right) ,  \label{3.21}
\end{equation}%
then for any isotonic linear functional $A:L\rightarrow \mathbb{R}$ with $%
A\left( h\right) >0$ we have the inequality%
\begin{equation}
0\leq A\left( hf^{2}\right) A\left( hg^{2}\right) -\left[ A\left( hfg\right) %
\right] ^{2}\leq \frac{1}{4}\left( M-m\right) ^{2}A^{2}\left( hg^{2}\right) .
\label{3.3}
\end{equation}%
The constant $\frac{1}{4}$ in (\ref{3.3}) is sharp.
\end{proposition}

\begin{proof}
We observe that 
\begin{equation*}
\left( Mg-f,f-mg\right) _{A,h}=A\left[ h\left( Mg-f\right) \left(
f-mg\right) \right] \geq 0.
\end{equation*}%
Applying Theorem \ref{t2.1} for $\left( \cdot ,\cdot \right) _{A,h}$ we get%
\begin{equation*}
0\leq \left( f,f\right) _{A,h}\left( g,g\right) _{A,h}-\left( f,g\right)
_{A,h}^{2}\leq \frac{1}{4}\left( M-m\right) ^{2}\left( g,g\right) _{A,h}^{2},
\end{equation*}%
which is clearly equivalent to (\ref{3.3}).
\end{proof}

The following corollary holds.

\begin{corollary}
\label{c3.2}Let $f,g\in F\left( T\right) $ such that $fg,$ $f^{2},g^{2}\in
F\left( T\right) .$ If $m,M$ are real numbers such that (\ref{3.21}) holds,
then%
\begin{equation}
0\leq A\left( f^{2}\right) A\left( g^{2}\right) -A^{2}\left( fg\right) \leq 
\frac{1}{4}\left( M-m\right) ^{2}A^{2}\left( g^{2}\right) .  \label{3.4}
\end{equation}%
The constant $\frac{1}{4}$ is sharp in (\ref{3.4}).
\end{corollary}

\begin{remark}
\label{r3.3}The condition (\ref{3.21}) may be replaced with the weaker
assumption%
\begin{equation}
\left( Mg-f,f-mg\right) _{A,h}\geq 0.  \label{3.5}
\end{equation}
\end{remark}

\section{Applications for Integrals}

Let $\left( \Omega ,\Sigma ,\mu \right) $ be a measure space consisting of a
set $\Omega ,$ a $\sigma -$algebra $\Sigma $ of subsets of $\Omega $ and a
countably additive and positive measure on $\Sigma $ with values in $\mathbb{%
R}\cup \left\{ \infty \right\} .$

Denote by $L_{\rho }^{2}\left( \Omega ,\mathbb{K}\right) $ the Hilbert space
of all $\mathbb{K}$-valued functions $f$ defined on $\Omega $ that are $%
2-\rho -$integrable on $\Omega ,$ i.e., $\int_{\Omega }\rho \left( t\right)
\left\vert f\left( s\right) \right\vert ^{2}d\mu \left( s\right) <\infty ,$
where $\rho :\Omega \rightarrow \lbrack 0,\infty )$ is a measurable function
on $\Omega .$

The following proposition contains a counterpart of the weighted
Cauchy-Buniakowsky-Schwarz's integral inequality.

\begin{proposition}
\label{p4.1}Let $A,a\in \mathbb{K}$ $\left( \mathbb{K}=\mathbb{C},\mathbb{R}%
\right) $ and $f,g\in L_{\rho }^{2}\left( \Omega ,\mathbb{K}\right) .$ If 
\begin{equation}
\int_{\Omega }\func{Re}\left[ \left( Ag\left( s\right) -f\left( s\right)
\right) \left( \overline{f\left( s\right) }-\overline{a}\text{ }\overline{g}%
\left( s\right) \right) \right] \rho \left( s\right) d\mu \left( s\right)
\geq 0  \label{4.1}
\end{equation}%
or, equivalently,%
\begin{equation*}
\int_{\Omega }\rho \left( s\right) \left\vert f\left( s\right) -\frac{a+A}{2}%
g\left( s\right) \right\vert ^{2}d\mu \left( s\right) \leq \frac{1}{4}%
\left\vert A-a\right\vert ^{2}\int_{\Omega }\rho \left( s\right) \left\vert
g\left( s\right) \right\vert ^{2}d\mu \left( s\right) ,
\end{equation*}%
then one has the inequality%
\begin{align}
0& \leq \int_{\Omega }\rho \left( s\right) \left\vert f\left( s\right)
\right\vert ^{2}d\mu \left( s\right) \int_{\Omega }\rho \left( s\right)
\left\vert g\left( s\right) \right\vert ^{2}d\mu \left( s\right) -\left\vert
\int_{\Omega }\rho \left( s\right) f\left( s\right) \overline{g\left(
s\right) }d\mu \left( s\right) \right\vert ^{2}  \label{4.2} \\
& \leq \frac{1}{4}\left\vert A-a\right\vert ^{2}\left( \int_{\Omega }\rho
\left( s\right) \left\vert g\left( s\right) \right\vert ^{2}d\mu \left(
s\right) \right) ^{2}.  \notag
\end{align}
\end{proposition}

\begin{proof}
Follows by Theorem \ref{t2.1} applied for the inner product $\left\langle
\cdot ,\cdot \right\rangle _{\rho }:=L_{\rho }^{2}\left( \Omega ,\mathbb{K}%
\right) \times L_{\rho }^{2}\left( \Omega ,\mathbb{K}\right) \rightarrow 
\mathbb{K}$,%
\begin{equation*}
\left\langle f,g\right\rangle _{\rho }:=\int_{\Omega }\rho \left( s\right)
f\left( s\right) \overline{g\left( s\right) }d\mu \left( s\right) .
\end{equation*}
\end{proof}

\begin{remark}
\label{r4.2}A sufficient condition for (\ref{4.1}) to hold is 
\begin{equation}
\func{Re}\left[ \left( Ag\left( s\right) -f\left( s\right) \right) \left( 
\overline{f\left( s\right) }-\overline{a}\text{ }\overline{g}\left( s\right)
\right) \right] \geq 0\text{ \ for \ }\mu -\text{a.e. }s\in \Omega .
\label{4.3}
\end{equation}
\end{remark}

In the particular case $\rho =1,$ we have the following counterpart of the
Cauchy-Buniakowsky-Schwarz inequality.

\begin{corollary}
\label{c4.3}Let $a,A\in \mathbb{K}$ $\left( \mathbb{K}=\mathbb{C},\mathbb{R}%
\right) $ and $f,g\in L_{\rho }^{2}\left( \Omega ,\mathbb{K}\right) .$ If%
\begin{equation}
\int_{\Omega }\func{Re}\left[ \left( Ag\left( s\right) -f\left( s\right)
\right) \left( \overline{f\left( s\right) }-\overline{a}\text{ }\overline{g}%
\left( s\right) \right) \right] d\mu \left( s\right) \geq 0,  \label{4.4}
\end{equation}%
or, equivalently%
\begin{equation*}
\int_{\Omega }\left\vert f\left( s\right) -\frac{a+A}{2}g\left( s\right)
\right\vert ^{2}d\mu \left( s\right) \leq \frac{1}{4}\left\vert
A-a\right\vert ^{2}\int_{\Omega }\left\vert g\left( s\right) \right\vert
^{2}d\mu \left( s\right) ,
\end{equation*}%
then one has the inequality%
\begin{eqnarray}
0 &\leq &\int_{\Omega }\left\vert f\left( s\right) \right\vert ^{2}d\mu
\left( s\right) \int_{\Omega }\left\vert g\left( s\right) \right\vert
^{2}d\mu \left( s\right) -\left\vert \int_{\Omega }f\left( s\right) 
\overline{g\left( s\right) }d\mu \left( s\right) \right\vert ^{2}
\label{4.5} \\
&\leq &\frac{1}{4}\left\vert A-a\right\vert ^{2}\left( \int_{\Omega
}\left\vert g\left( s\right) \right\vert ^{2}d\mu \left( s\right) \right)
^{2}.  \notag
\end{eqnarray}
\end{corollary}

\begin{remark}
\label{r4.4}If $\mathbb{K}=\mathbb{R},$ then a sufficient condition for
either (\ref{4.1}) or (\ref{4.4}) \ to hold is 
\begin{equation}
ag\left( s\right) \leq f\left( s\right) \leq Ag\left( s\right) \text{ \ for
\ }\mu -\text{a.e. }s\in \Omega ,  \label{4.6}
\end{equation}%
where, in this case, $a,A\in \mathbb{R}$ with $A>a.$
\end{remark}

\section{Applications for Sequences}

For a given sequence $\left( w_{i}\right) _{i\in \mathbb{N}}$ of nonnegative
real numbers, consider the Hilbert space $\ell _{w}^{2}\left( \mathbb{K}%
\right) ,$ $\left( \mathbb{K}=\mathbb{C},\mathbb{R}\right) ,$ where 
\begin{equation}
\ell _{w}^{2}\left( \mathbb{K}\right) :=\left\{ \overline{x}=\left(
x_{i}\right) _{i\in \mathbb{N}}\subset \mathbb{K}\left\vert
\sum_{i=0}^{\infty }w_{i}\left\vert x_{i}\right\vert ^{2}<\infty \right.
\right\} .  \label{5.1}
\end{equation}

The following proposition that provides a counterpart of the weighted
Cauchy-Bunyakowsky-Schwarz inequality for complex numbers holds.

\begin{proposition}
\label{p5.1}Let $a,A\in \mathbb{K}$ and $\overline{x},\overline{y}\in \ell
_{w}^{2}\left( \mathbb{K}\right) .$ If%
\begin{equation}
\sum_{i=0}^{\infty }w_{i}\func{Re}\left[ \left( Ay_{i}-x_{i}\right) \left( 
\overline{x_{i}}-\overline{a}\text{ }\overline{y_{i}}\right) \right] \geq 0,
\label{5.2}
\end{equation}%
then one has the inequality%
\begin{equation}
0\leq \sum_{i=0}^{\infty }w_{i}\left\vert x_{i}\right\vert
^{2}\sum_{i=0}^{\infty }w_{i}\left\vert y_{i}\right\vert ^{2}-\left\vert
\sum_{i=0}^{\infty }w_{i}x_{i}\overline{y_{i}}\right\vert ^{2}\leq \frac{1}{4%
}\left\vert A-a\right\vert ^{2}\left( \sum_{i=0}^{\infty }w_{i}\left\vert
y_{i}\right\vert ^{2}\right) ^{2}.  \label{5.3}
\end{equation}%
The constant $\frac{1}{4}$ is sharp.
\end{proposition}

\begin{proof}
Follows by Theorem \ref{t2.1} applied for the inner product $\left\langle
\cdot ,\cdot \right\rangle _{w}:\ell _{w}^{2}\left( \mathbb{K}\right) \times
\ell _{w}^{2}\left( \mathbb{K}\right) \rightarrow \mathbb{K}$,%
\begin{equation*}
\left\langle \overline{x},\overline{y}\right\rangle _{w}:=\sum_{i=0}^{\infty
}w_{i}x_{i}\overline{y_{i}}.
\end{equation*}
\end{proof}

\begin{remark}
\label{r5.2}A sufficient condition for (\ref{5.2}) to hold is%
\begin{equation}
\func{Re}\left[ \left( Ay_{i}-x_{i}\right) \left( \overline{x_{i}}-\overline{%
a}\overline{y_{i}}\right) \right] \geq 0\text{ \ for all }i\in \mathbb{N}.
\label{5.4}
\end{equation}
\end{remark}

In the particular case $w_{i}=1,$ $i\in \mathbb{N}$, we have the following
counterpart of the Cauchy-Bunyakowsky-Schwarz inequality.

\begin{corollary}
\label{c5.3}Let $a,A\in \mathbb{K}$ $\left( \mathbb{K}=\mathbb{C},\mathbb{R}%
\right) $ and $\overline{x},\overline{y}\in \ell ^{2}\left( \mathbb{K}%
\right) .$ If%
\begin{equation}
\sum_{i=0}^{\infty }\func{Re}\left[ \left( Ay_{i}-x_{i}\right) \left( 
\overline{x_{i}}-\overline{a}\overline{y_{i}}\right) \right] \geq 0,
\label{5.5}
\end{equation}%
then one has the inequality%
\begin{equation}
0\leq \sum_{i=0}^{\infty }\left\vert x_{i}\right\vert ^{2}\sum_{i=0}^{\infty
}\left\vert y_{i}\right\vert ^{2}-\left\vert \sum_{i=0}^{\infty }x_{i}%
\overline{y_{i}}\right\vert ^{2}\leq \frac{1}{4}\left\vert A-a\right\vert
^{2}\left( \sum_{i=0}^{\infty }\left\vert y_{i}\right\vert ^{2}\right) ^{2}.
\label{5.6}
\end{equation}
\end{corollary}

\begin{remark}
\label{r5.4}If $\mathbb{K}=\mathbb{R}$, then a sufficient condition for
either (\ref{5.2}) or (\ref{5.5}) to hold is%
\begin{equation}
ay_{i}\leq x_{i}\leq Ay_{i}\text{ \ for each \ }i\in \mathbb{N},  \label{5.7}
\end{equation}%
with $A>a.$
\end{remark}

\end{document}